\documentclass[12pt]{article}

\usepackage{moreverb}

\usepackage[colorlinks,bookmarksopen,bookmarksnumbered,citecolor=red,urlcolor=red]{hyperref}
\usepackage{amssymb,amsmath,amsthm}
\usepackage{color}
\usepackage{graphicx}
\usepackage{subfigure}
\usepackage{epstopdf}
\usepackage{graphics,epsfig,cite}

\newcommand\BibTeX{{\rmfamily B\kern-.05em \textsc{i\kern-.025em b}\kern-.08em
T\kern-.1667em\lower.7ex\hbox{E}\kern-.125emX}}

\begin{document}


\title{Real-Time Non-Linear Receding Horizon Control Methodology for Estimation of Time-Varying Parameters }

\author{Fei Sun*, Kamran Turkoglu \\ (\emph{submitted and under review in Optimal Control Applications and Methods})}


\maketitle
\begin{abstract}
In control and engineering community, models generally contain a number of parameters which are unknown or roughly known. A complete knowledge of these parameters is critical to describe and analyze the dynamics of the system. This paper develops a novel approach of estimating unknown time-varying parameters of nonlinear systems based on real-time nonlinear receding horizon control strategy. For this purpose, the problem of estimation is converted into optimization problem which can be solved by backwards sweep Riccati method. Methodology is implemented on two nonlinear examples without providing any reference of the parameters. Results successfully demonstrate the fact that proposed algorithm is able to estimate unknown time-varying parameters of nonlinear systems effectively.
\end{abstract}



\vspace{-6pt}

\section{Introduction}
\vspace{-2pt}
Parameter estimation is a widely used process in modeling of dynamic systems, where systems exhibit unknown but desired dynamics. In real world applications, many dynamical systems, including systems which belong to the fields of biology, chemistry, physics, engineering and many more \cite{adetola2009,ho2010,wang2012,lin2014,bellsky2014,ding2014,guay2014acc,tutsoy2015},  may (and mostly do) have partial or all unknown parameters. Knowledge about the parameters of the underlying dynamics is the prerequisite to analyze, control, and predict their characteristic behaviors. Hence, this topic has drawn great attention in various areas due to its theoretical and practical significance.

While estimating unknown parameters in dynamic systems could be a daunting task, study of parameter estimation, especially in nonlinear systems that exhibit time-varying behavior, remains an open field of research. In this area, an important effort is the application of adaptive control and synchronization methodologies by designing an adaptive controller to estimate the uncertainties as well as minimize the synchronization error \cite{wu1993,li2004,huang2006pre,li2007,sun2009,sun2012}. However, adaptive control based methodologies highly depend on the assumption that parameters are constant or slowly time-varying. In practice, the time-varying feature of system parameters can have significant impact on the performance of the underlying systems \cite{dochain2003,bravo2006,butt2013,wang2014nd}, and canot be underestimated.

On the other hand, especially with the improvements in microprocessors and increased computational capabilities, Receding Horizon Control (RHC) based methodologies gained significant momentum in the last two decades \cite{bryson1975,thomas1975,chen1982,kwon1983,mayne1990,mayne1990b,he2008}. RHC is a formulation which is widely used to obtain an optimal feedback control law by minimizing the desired performance index for a given finite horizon. The performance index of a receding-horizon control problem has a moving initial time and a moving terminal time, and the time interval of the performance index is finite. Since the time interval of the performance index is finite, the optimal feedback law can be determined even for a system that is not stabilizable. One advantage of receding horizon optimal control formulation is that it can deal with a broader class of control objectives other than asymptotic stabilization. Moreover, nonlinear receding horizon control (NRHC) has made a significant impact on industrial controls and is being increasingly applied in process controls \cite{he2008,mayne1990,mayne1990b}. Various advantages are known for NRHC, including the ability to handle time-varying and nonlinear systems, input/output constraints and plant/parameter related uncertainty.

There are many existing valuable works related to the use of different methodologies in nonlinear systems for parameter estimation procedures. One of the conducted studies concentrates on parameter identification using differential evolution algorithms \cite{ho2010}, while another study provides a bio-geography based parameter estimation for nonlinear systems \cite{wang2012}. Other important studies include Kalman filtering based parameter estimation \cite{bellsky2014}, adaptive estimation with invariant manifolds \cite{guay2014acc}, estimation through identical or non-identical structures \cite{sun2012}, chaotic ant-swarm models \cite{li2010pre,li2006}, time series approach \cite{yu2007pre,parlitz1996}, via particle swarm optimization \cite{he2007}, recursive identification via incremental estimation \cite{zhou1996}.

The main goal in this approach remains as the parameter estimation of nonlinear time varying systems, but the fundamental difference of this study with respect to the existing literature is that in this paper it is desired to demonstrate the connection (and applicability) of real-time nonlinear receding horizon control algorithms as an effective parameter estimation procedure. Here, we provide a framework which is able to cope with nonlinear time-varying systems. For this purpose, the estimation procedure is reduced to a family of finite horizon optimization problems. In addition, a non-iterative optimization algorithm is employed to avoid high computational complexity and large data storage. With this formulation, we also provide the closed-loop asymptotic stability guarantees of the synchronization error system. Obtained results from two example problems demonstrate the applicability of the real-time NRHC as a parameter estimation routine on the nonlinear systems with unknown time-varying parameters. This approach also serves as a minimization routine the synchronization errors, as well.

The paper is organized as the following: Problem formulation is defined in Section-\ref{sec:prob_formulation}.  In Section-\ref{sec:backward_sweep}, backward-sweep method is iterated and stability proof is provided. In Section-\ref{sec:results}, we discuss the findings of such approach, and with the conclusions section, we finalize the paper.


\section{Problem formulation}\label{sec:prob_formulation}
\vspace{-2pt}

For this study, we consider the following nonlinear system with unknown time-varying parameters:
\begin{equation}\label{eq:sys_eom}
\dot{x}=Ax+f(x)+D(x)\Theta(t),\\
\end{equation}
where $x\in R^{n}$ is the state vector, $A\in R^{n\times n}$ is the the linear coefficient
matrix and $f(x):R^{n}\rightarrow R^{n}$ is the nonlinear part of system in \eqref{eq:sys_eom}. $D(x):R^{n}\rightarrow R^{n\times p}$ is a known function vector and $\hat{\Theta} \in R^{p}$ denotes the unknown time-varying parameters.

In order to estimate the unknown parameters, the corresponding controlled system is given by

\begin{equation}\label{eq:sys_eom_response}
\dot{y}=By+f(y)+D(y)\hat{\Theta}(t),\\
\end{equation}
where $y\in R^{n}$ is the state vector and $\hat{\Theta}$ represents the estimated parameter. Here, function $f(\cdot)$ and $D(\cdot)$ satisfy the global Lipschitz condition. Therefore there exist positive constants $\alpha$ and $\beta$ such that
\begin{equation*}
\begin{split}
\|f(y)-f(x)\|\le&\beta_1\|y-x\|,\\
\|D(y)-D(x)\|\le&\beta_2\|y-x\|.
\end{split}
\end{equation*}

The system in Eq.\eqref{eq:sys_eom} is considered as the drive system and system in Eq.\eqref{eq:sys_eom_response} is considered as the response system. The synchronization error between the drive and the response system becomes an important part of this analysis, and in this study it is defined as
$$e(t) = y(t) -x(t).$$
With this definition, it is possible to define the \emph{synchronization error system} as
\begin{equation}
\dot{e}(t)=Ae(t)+B(f(y(t))-f(x(t)))+D(y)\hat{\Theta}(t)-D(x)\Theta(t)
\end{equation}
where the estimation error is denoted by
$$\tilde{\Theta}(t)=\hat{\Theta}(t)-\Theta(t).$$

From this point on, for the estimation of unknown time-varying parameters, we propose and establish a finite horizon cost function which is associated with the synchronization error and the estimated parameter, as:

\begin{equation}\label{eq:sync_cost_func}
J=\frac{1}{2}\int_t^{t+T}[e^TQe+\hat{\Theta}^TR\hat{\Theta}]{\rm d}\tau,
\end{equation}
where weighting matrices $Q>0$ and $R>0$. With the construction of the synchronization cost function, as given in Eq.\eqref{eq:sync_cost_func}, the estimation problem is converted into a parameter optimization procedure, where the unknown parameter(s) can be estimated through this process.
For this purpose, we utilize the powerful nature of \emph{real-time nonlinear receding horizon control algorithm} to minimize associated synchronization the cost function. In this context, the performance index evaluates the performance from the present time $t$ to the finite future $t+T$, and is minimized for each time segment $t$ starting from $y(t)$. With this structure, it is possible to convert the present receding horizon control problem into a family of finite horizon optimal control problems on the artificial $\tau$ axis parameterized by time $t$.

It is well known from literature that first order necessary conditions of optimality are obtained from the two-point boundary value problem (TPBVP) \cite{bryson1975}, as

\begin{equation}\label{eq:tpbvp}
\begin{split}
y_\tau^* = H_{\lambda}^T , \hskip 20pt y^* = &y(t), \hskip 20pt \lambda^*(\tau,t) = -H_y ^T \\
 \lambda^*(T,t) =&~ 0, \hskip 20pt  H_{\hat{\Theta}} = 0
\end{split}
\end{equation}
where H is the Hamiltonian and is defined as
\begin{equation}\label{eq:tpbvp_ham}
\begin{split}
H &= L + \lambda^{*T}\dot{y}\\
&=\frac{1}{2}(e^TQe+\hat{\Theta}^TR\hat{\Theta})+ \lambda^{*T}[Ay+f(y)+D(y)\hat{\Theta}(t)].
\end{split}
\end{equation}
In Eqs.\eqref{eq:tpbvp}-\eqref{eq:tpbvp_ham}, $(~~)^*$ denotes a variable in the optimal control problem so as to distinguish it from its correspondence in the original problem.

Using this formulation, the estimated parameter can be calculated as
\begin{equation}\label{theta}
\hat{\Theta}(t) = arg\{H_{\hat{\Theta}}[y(t),\lambda(t),\hat{\Theta}(t),x(t)] = 0\}
\end{equation}

In this methodology, the TPBVP is to be regarded as a nonlinear equation with respect to the costate at $\tau=0$ as
\begin{equation}
F(\lambda(t),y(t),T,t)=\lambda^*(T,t)=0.
\end{equation}

Since the nonlinear equation $F(\lambda(t),y(t),T,t)$ has to be satisfied at any time $t$, $\frac{{\rm d}F}{{\rm d}t}=0$ holds along the trajectory of the closed-loop system of the receding horizon control. If $T$ is a smooth function of time $t$, the solution of $F(\lambda(t),y(t),T(t),t)$ can be tracked with respect to time. In this formulation, the ordinary differential equation of $\lambda(t)$ can be solved numerically, in real-time, without any need of an iterative optimization routine. However, numerical errors associated with the solution may accumulate as the integration proceeds in practice, and therefore some correction techniques are required to correct such errors in the solution. To address this problem, a stabilized continuation method \cite{kabamba1987,ohtsuka1994,ohtsuka1997,ohtsuka1998} is used. According to this method, it is possible to rewrite the statement as
\begin{equation}\label{cm}
\frac{{\rm d}F}{{\rm d}t}=-A_sF,
\end{equation}
where $A_s$ denotes a stable matrix to make the solution converge to zero exponentially, such that
\begin{equation}
F=F_0e^{-A_st}\rightarrow 0 \,(t\rightarrow \infty),
\end{equation}
where $F_0$ denotes the initial value of $F$. And thus the numerical errors is attenuated through the integration process.

\section{Backward-sweep Method:}\label{sec:backward_sweep}

In order to integrate the differential equation of $\lambda(t)$ in real time, the partial differentation of Eqs.\eqref{eq:tpbvp} with respect to time $t$ and $\tau$, converts the problem to the following linear differential equation:
\begin{equation}\label{pd}
\frac {\partial}{\partial \tau}\begin{bmatrix}y^*_t-y^*_{\tau}\\
\lambda^*_t-\lambda^*_{\tau} \end{bmatrix}=\begin{bmatrix}G&-L\\
-K&-G^T \end{bmatrix}\begin{bmatrix}y^*_t-y^*_{\tau}\\
\lambda^*_t-\lambda^*_{\tau} \end{bmatrix}
\end{equation}
where $G=f_y-f_{\hat{\Theta}}H^{-1}_{\hat{\Theta}\hat{\Theta}}H_{\hat{\Theta} y}$, $L=f_{\hat{\Theta}}H^{-1}_{\hat{\Theta}\hat{\Theta}}f_{\hat{\Theta}}^T$, $K=H_{yy}-H_{y\hat{\Theta}}H^{-1}_{\hat{\Theta}\hat{\Theta}}H_{\hat{\Theta} y}$. Since $x_t$ and $x_{\tau}$ are canceled in Eq.\eqref{pd}, data storage is reduced.

The derivative of the nonlinear function $F$ with respect to time is given by
\begin{equation}\label{df}
\frac{{\rm d}F}{{\rm d}t}=\lambda^*_t(T,t)+\lambda^*_{\tau}(T,t)\frac{{\rm d}T}{{\rm d}t}.
\end{equation}

In order to reduce the computational cost without resourcing any approximation technique, the backward-sweep method is implemented. This provides the relationship between the costate and other variables is expressed as the followings:
\begin{equation}\label{relation}
\lambda^*_t-\lambda^*_\tau=S(\tau,t)(y^*_t-y^*_\tau)+c(\tau,t),
\end{equation}
where
\begin{equation}\label{sc}
\begin{split}
S_{\tau}&=-G^TS-SG+SLS-K, \\
c_{\tau}&=-(G^T-SL)c,
\end{split}
\end{equation}
and
\begin{equation}\label{sct}
\begin{split}
S(T,t)&=0,\\
c(T,t)&=H_y^T\mid_{\tau=T}(1+\frac{{\rm d}T}{{\rm d}t})-A_sF.
\end{split}
\end{equation}
to be satisfied. With this analogy, the differential equation of $\lambda(t)$ to be integrated in real time is obtained as:
\begin{align}\label{dl}
\frac{{\rm d}\lambda(t)}{{\rm d}t}=-H_y^T+c(0,t).
\end{align}

Here, at each time $t$, the Euler-Lagrange equations (Eqs.\eqref{eq:tpbvp}) are integrated forward along the $\tau$ axis, and then Eqs.\eqref{sc} are integrated backward with terminal conditions provided in Eq.\eqref{sct}. Next, the differential equation of $\lambda(t)$ is integrated for one step along the $t$ axis so as to estimate the unknown parameters from Eq.\eqref{theta}. If the matrix $H_{\hat{\Theta}\hat{\Theta}}$ is nonsingular, the algorithm is executable regardless of controllability or stabilizability or the system.

\subsection{Convergence and Stability Analysis of NRHC Estimation Problem}

\textbf{Theorem-1:} If $\exists$  $V(e,\hat{\Theta})$ such that $\dot{V}(e,\hat{\Theta})<0$ $\forall$ $V(0)=0$, $V>0$ and $e\ne 0$ , then the closed loop dynamics defined for the drive system $\dot{x}=Ax+f(x)+D(x)\Theta(t)$, the response system $\dot{y}=By+f(y)+D(y)\hat{\Theta}(t)$ and the error system $\dot{e}(t)=Ae(t)+B(f(y(t))-f(x(t)))+D(y)\hat{\Theta}(t)-D(x)\Theta(t)$, using the nonlinear receding horizon control methodology,  is asymptotically stable (i.e. there exists a stability ball $Z(0,z)) \in R^n$ $\forall$ $z>0$, such that for any given initial value, the synchronization error ($e$) is stable and converges ($e\rightarrow 0$) as $t \rightarrow \infty$.)

\begin{proof}
Let's take into account the following Lyapunov function with respect to $e$ and $\hat{\Theta}$ in the form of
\begin{equation}
V=\frac{1}{2}\int_t^{t+T}(e^TQe+\hat{\Theta}^T R\hat{\Theta})d\tau.
\end{equation}
where clearly we have $V(0)=0$ and $V>0$ $\forall$ $e\neq 0$.

Here, it is possible to partition the candidate Lyapunov function as
\begin{equation}
V(e^*(t))=\frac{1}{2}\int_t^{t+\Delta t}(e^T(\tau)Qe(\tau)+\hat{\Theta}^T(\tau) R\hat{\Theta}(\tau))d\tau+\frac{1}{2}\int_{t+\Delta t}^{t+T}(e^T(\tau)Qe(\tau)+\hat{\Theta}^T(\tau) R\hat{\Theta}(\tau))d\tau,
\end{equation}
where $e(\tau)=e(\tau;e^*(t),t)$ and $\hat{\Theta}(\tau)=\hat{\Theta}(\tau;\hat{\Theta}^*(t),t)$ are optimal values for system, at given time instant $t$. Now, consider a parameter estimator defined as follows:
\begin{equation}
\hat{\Theta}(\tau)=
\begin{cases}
\hat{\Theta}(\tau;\hat{\Theta}^*(t),t)& \text{for $\tau \in [t+\Delta t,t+T]$},\\
0& \text{for $\tau \in (t+T,t+T+\Delta t)$}.
\end{cases}
\end{equation}
Let $\bar{e}(\cdot)=\bar{e}(\cdot:e(t+\Delta t),t+\Delta t):[t+\Delta t, t+T+\Delta t]\in R^{n}$ and denote the corresponding error trajectory with initial condition $\bar{e}(t+\Delta t)=\bar{e}(t+\Delta t;e(t),t)$. It is clear that
\begin{equation}
\bar{e}(\tau)=
\begin{cases}
e(\tau;e^*(t),t)& \text{for $\tau \in [t+\Delta t,t+T]$},\\
0&  \text{for $\tau \in (t+T,t+T+\Delta t)$}.
\end{cases}
\end{equation}
because $\bar{e}(t+T;e(t+\Delta t),t+\Delta t)=e(t+T;e^*(t),t)=0$ and $\bar{\hat{\Theta}}=0$ for $\tau>t+T$. Since $\bar{\hat{\Theta}}$ is not necessarily optimal for trajectory of error $e(t+\Delta t)$ from time $t+\Delta t$ \cite{mayne1990}, it follows that
\begin{equation}
\begin{split}
V(e^*(t))&=\frac{1}{2}\int_t^{t+\Delta t}(e^T(\tau)Qe(\tau)+\hat{\Theta}^T(\tau) R\hat{\Theta}(\tau))d\tau+V(e(t+\Delta t), t+\Delta t;\bar{\hat{\Theta}})\\
&\geqslant \frac{1}{2}\int_t^{t+\Delta t}(e^T(\tau)Qe(\tau)+\hat{\Theta}^T(\tau) R\hat{\Theta}(\tau))d\tau+V(e(t+\Delta t))
\end{split}
\end{equation}
so that
\begin{equation}\label{dv}
V(e(t+\Delta t))-V(e^*(t)) \leq -\frac{1}{2}\int_t^{t+\Delta t}(e^T(\tau)Qe(\tau)+\hat{\Theta}^T(\tau) R\hat{\Theta}(\tau))d\tau.
\end{equation}

Since $V$ is continuously differentiable, from the Mean Value Theorem \cite{mayne1990} we have
\begin{equation}
\frac{V(e(t+\Delta t))-V(e^*(t))}{\Delta t}=\triangledown_eV(e^*(t))+\rho(\Delta t)(e(t+\Delta t)-e^*(t))\frac{e(t+\Delta t)-e^*(t)}{\Delta t}
\end{equation}
where $\rho(\Delta t)\in (0,1)$. Since
\begin{equation}
e^*(t)=e(t;e^*(t),t)=e(t), \, \hat{\Theta}^*(t)=\hat{\Theta}(t;e^*(t),t)=\hat{\Theta}(t)
\end{equation}
and $\hat{\Theta}$ is continuous at $t$, we have that
\begin{equation}
\lim_{\Delta t\rightarrow 0_+}\frac{e(t+\Delta t)-e^*(t)}{\Delta t}=f(y(t),\hat{\Theta}(t),x(t))=f(y^*(t),\hat{\Theta}^*(t)).
\end{equation}
Since $\triangledown_e V$ and $e$ are both continuous, we get
\begin{equation}
\lim_{\Delta t\rightarrow 0_+}\frac{V(e(t+\Delta t))-V(e^*(t))}{\Delta t}=\triangledown_eV(e^*(t))f(y^*(t),\hat{\Theta}^*(t),x(t)).
\end{equation}
Based on the above derivation, we also have
\begin{equation}
\begin{split}
\lim_{\Delta t\rightarrow 0_+}&-\frac{1}{2\Delta t}\int_t^{t+\Delta t}(e^T(\tau)Qe(\tau)+\hat{\Theta}^T(\tau) R\hat{\Theta}(\tau))d\tau\\
&\leq \lim_{\Delta t\rightarrow 0_+}-\frac{1}{2\Delta t}\int_t^{t+\Delta t}(e^T(\tau)Qe(\tau))d\tau\\
&=-\frac{1}{2}(e^T(\tau)Qe(\tau)).
\end{split}
\end{equation}
Dividing both sides of \eqref{dv} by $\Delta t >0$ and taking the limit as $\Delta t \rightarrow 0_+$ yields to
\begin{equation}
\triangledown_eV(e^*(t))f(y^*(t),\hat{\Theta}^*(t),y(t))\leq -\frac{1}{2}(e^T(\tau)Qe(\tau))\leq 0.
\end{equation}
Therefore,
\begin{equation}
\frac{d V(e^*(t))}{dt}=\triangledown_eV(e^*(t))f(y^*(t),\hat{\Theta}^*(t),x(t))\leq 0.
\end{equation}

Obviously, $\dot{V}(e(t)) \le 0$ and therefore $V(e(t))$ (namely, $J(e(t))$) is non-creasing.

It is clear that $\int_t^{t+T}(e^TQe+\hat{\Theta}^T R\hat{\Theta})d\tau \rightarrow 0$ as $t\rightarrow \infty$, which indicates the asymptotic stability of the synchronization error.
\end{proof}
With this result, we show that the non-increasing monotonicity of the cost function is a sufficient condition for the stability and the unknown time-varying parameters can be estimated through this process, without any given reference trajectory for the parameters.

\section{Simulation Results}\label{sec:results}

In this section, we use two benchmark examples to verify the effectiveness of the proposed method.

\textbf{Example-1:} We first consider the Lorenz system \cite{lorenz1963} with time-varying parameters.The objective is to estimate the time evolution of unknown parameters using the proposed nonlinear receding horizon algorithm and show the applicability of proposed algorithm on a nonlinear dynamics (such as modified Lorenz system).

Here, modified Lorenz system with time-varying parameters is given by
\begin{equation}\label{system}
\begin{split}
&\frac{{\rm d}x(t)}{{\rm d}t}=\begin{bmatrix}\frac{10\sin(t)}{t+1}(x_2-x_1) \\ 28x_1-x_1x_3-x_2\\x_1x_2-\frac{8}{3}x_3\end{bmatrix},\\
&\frac{{\rm d}y(t)}{{\rm d}t}=\begin{bmatrix}\hat{\theta}_1(y_2-y_1) \\ 28y_1-y_1y_3-y_2\\y_1y_2-\hat{\theta}_2y_3\end{bmatrix}.
\end{split}
\end{equation}
For this specific problem, the performance index is chosen as follows:
\begin{equation}
J=\frac{1}{2}\int_t^{t+T}\{[y(t^{\prime})-x(t^{\prime})]^TQ[y(t^{\prime})-x(t^{\prime})]+\hat{\Theta}(t^{\prime})^TR\hat{\Theta}(t^{\prime})\}{\rm d}t^{\prime},
\end{equation}
where $Q=\text{diag}(1,0.5,0.1)$ and $R=\text{diag}(0.02,0.02)$.

The smooth function of time interval $T$, taken from \cite{ohtsuka1997}, in the performance index is given by
\begin{equation}
T(t)=T_f(1-e^{-\alpha t}),
\end{equation}
where $T_f=0.1$ and $\alpha=0.01$. The stable matrix is chosen as $A_s=50I$.

The initial states are given by
\begin{equation}\label{initial1}
\begin{bmatrix}x_1(0) \\ x_2(0)\\x_3(0)\end{bmatrix}=\begin{bmatrix}-3 \\ -3\\15\end{bmatrix},
\begin{bmatrix}y_1(0) \\ y_2(0)\\y_3(0)\end{bmatrix}=\begin{bmatrix}-6 \\ -6\\22\end{bmatrix}.
\end{equation}

\begin{figure*}\label{fig:ex1_traj}
 \centering
 \vspace*{-5.cm}
    \includegraphics[width=12cm]{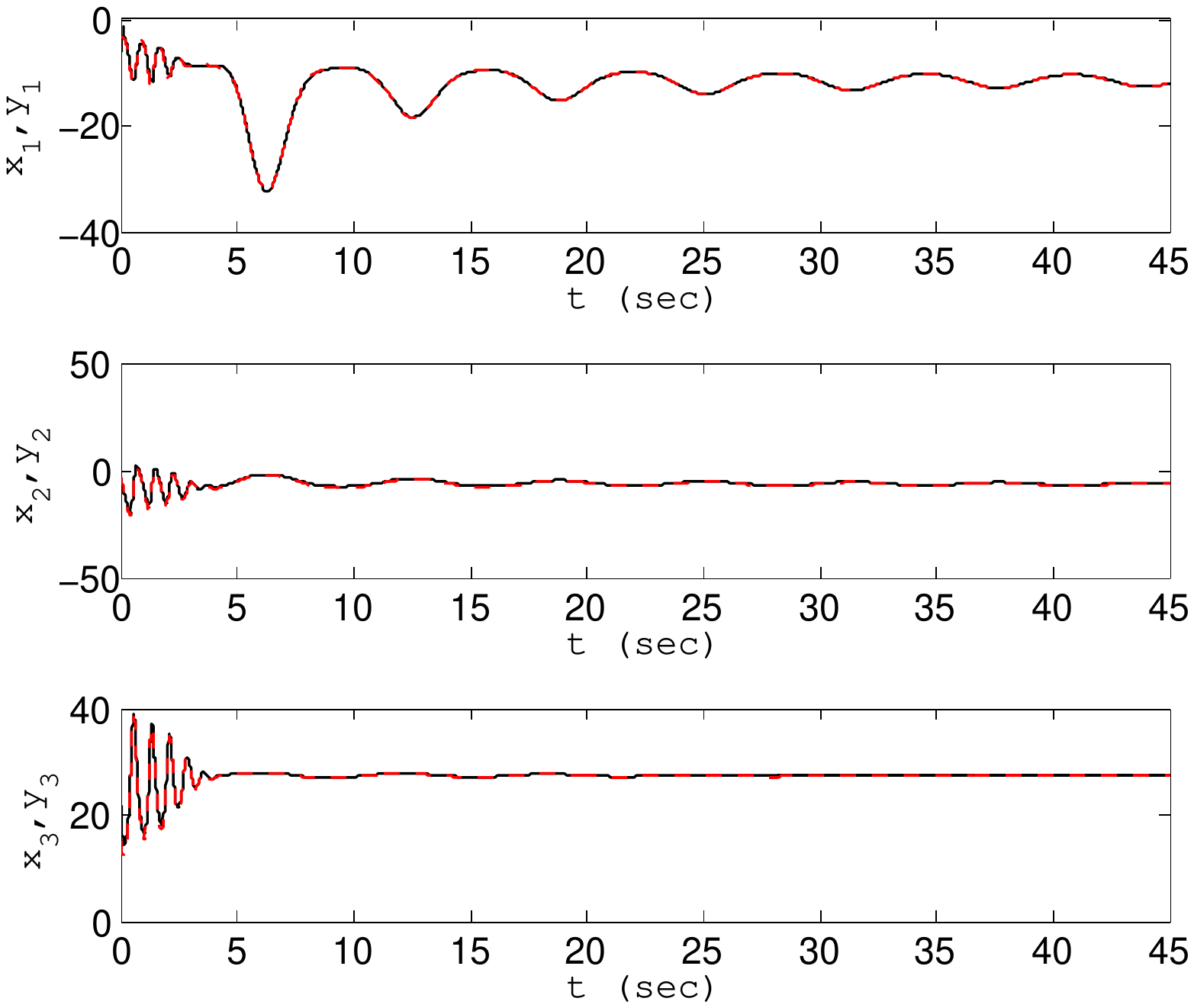}
    \vspace*{-5.cm}
    \caption{ \label{fig1} (Color online) The trajectories of states with time invariant parameters where solid
line denotes the trajectory of response system and dash line denotes the trajectory of drive system.}
 \end{figure*}

\begin{figure}
 \centering
 \vspace*{-1.cm}
    \includegraphics[width=14cm]{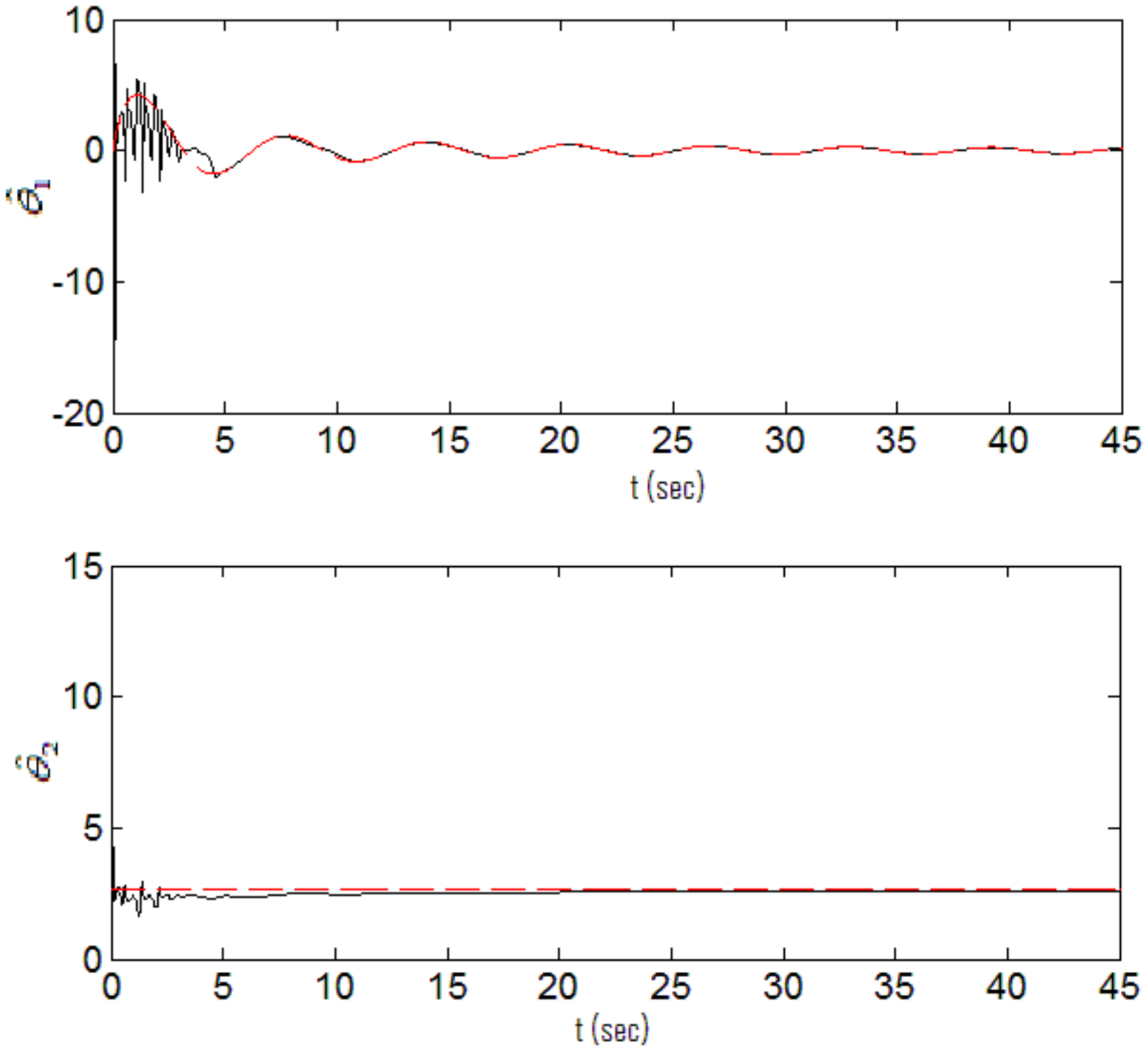}
    \vspace*{-1.cm}
    \caption{ \label{fig2} (Color online) The trajectories of estimated constant parameters where solid line
denotes the trajectory of estimated parameters and dash line denotes the trajectory of their
true values.}
 \end{figure}

\begin{figure*}
 \centering
    \includegraphics[width=14cm]{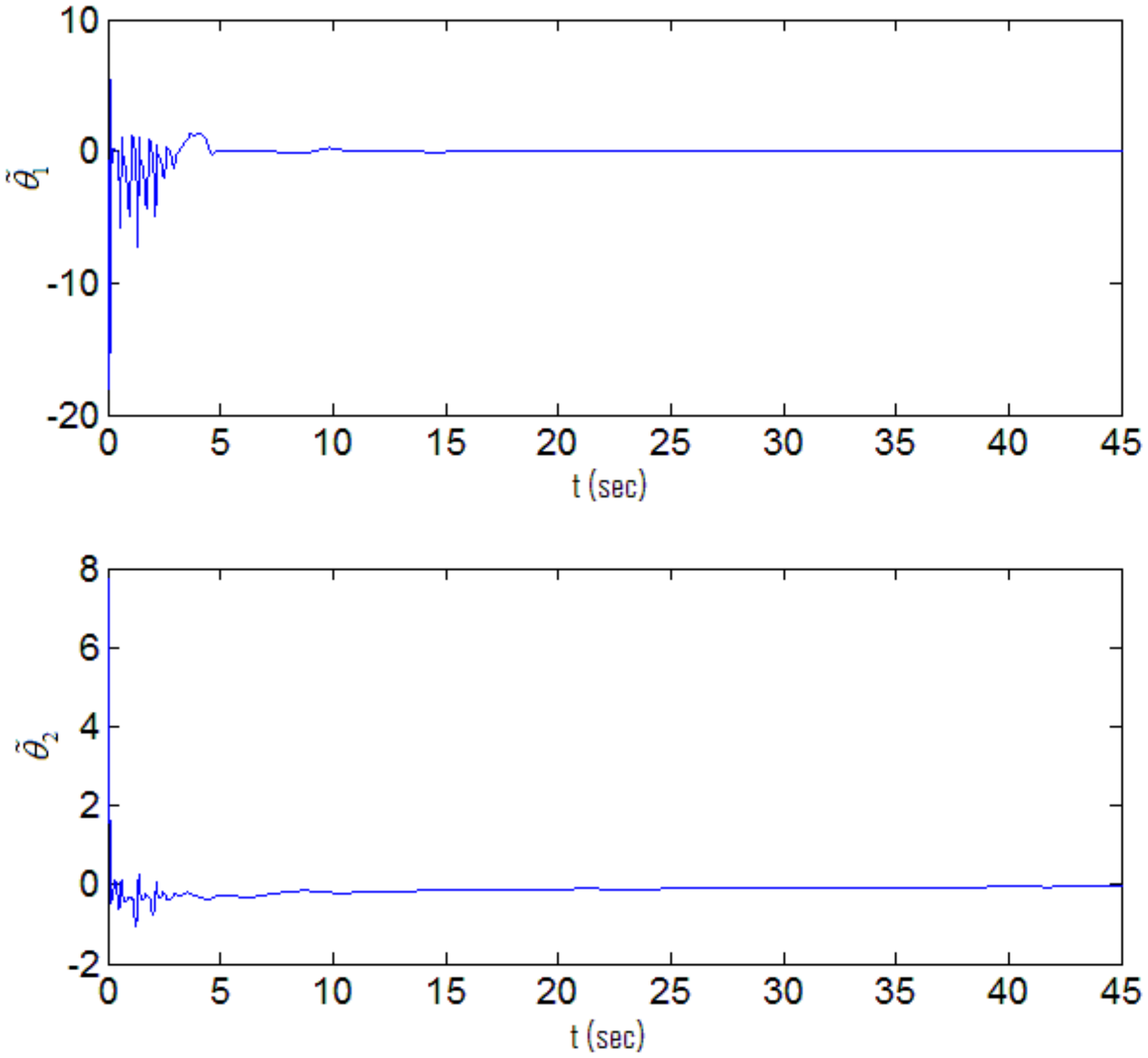}
    \vspace*{-1.cm}
    \caption{ \label{fig3} (Color online)The trajectories of estimated errors.}
 \end{figure*}

The simulation is implemented in MATLAB, where the time step on the t axis
is $0.01s$ and the time step on the artificial $\tau$ axis is $0.005s$. Fig.\ref{fig1} shows the trajectories
of drive-response systems in Eqs.\eqref{system} with initial conditions given in Eqs.\eqref{initial1}. The
estimated parameters $\hat{\theta}_1$ and $\hat{\theta}_2$ are presented in Fig.\ref{fig2}, where Fig.\ref{fig3} depicts the estimated error which clearly show the estimated parameters converge to their true values by using the NRHC
method.


It is possible to see from the Figs.\ref{fig1}-\ref{fig3} that the convergence rate of the estimated parameter is less than 5 seconds, and is relatively fast. In a real-life application, this methodology (and the associated estimation results) could be used to provide (feedback) information regarding the unknown parameter after a grace period (say 8-10sec) and/or the desired level of convergence is achieved.

\textbf{Example-2:} In the following, we consider a nonlinear system taken from \cite{guay2014acc}:
\begin{equation}\label{system2}
\begin{split}
&\frac{{\rm d}x(t)}{{\rm d}t}=\begin{bmatrix}-x_2^2-2x_1+\theta_1(t) \\ -2x_2+\theta_2(t)x_1\end{bmatrix},\\
&\frac{{\rm d}y(t)}{{\rm d}t}=\begin{bmatrix}-y_2^2-2y_1+\hat{\theta}_1(t) \\ -2y_2+\hat{\theta}_2(t)y_1.\end{bmatrix},
\end{split}
\end{equation}
where the time-varying parameters are defined as
\begin{equation}
\begin{split}
&\theta_1(t)=\begin{cases}2+\sin(t)& 0\le t\le6\pi\\2-\sin(\frac{\pi}{3})+\sin(2t+\frac{\pi}{3})& 6\pi\le t\le14\pi+\frac{\pi}{12}\\2-\sin(\frac{\pi}{3})+\sin(\frac{\pi}{2})& t\ge14\pi+\frac{\pi}{12}\end{cases}\\
&\theta_2(t)=3\cos(0.1\pi t).
\end{split}
\end{equation}
and the initial states are given by
\begin{equation}\label{initial2}
\begin{bmatrix}x_1(0) \\ x_2(0)\end{bmatrix}=\begin{bmatrix}0 \\ 0\end{bmatrix},
\begin{bmatrix}y_1(0) \\ y_2(0)\end{bmatrix}=\begin{bmatrix}1 \\ 2\end{bmatrix}.
\end{equation}

For this specific problem, the performance index is chosen as follows:
\begin{equation}
J=\frac{1}{2}\int_t^{t+T}\{[y(t^{\prime})-x(t^{\prime})]^TQ[y(t^{\prime})-x(t^{\prime})]+\hat{\Theta}(t^{\prime})^TR\hat{\Theta}(t^{\prime})\}{\rm d}t^{\prime},
\end{equation}
where $Q=\text{diag}(100,110)$ and $R=\text{diag}(0.1,0.2)$. The stable matrix is designed as $A_s=10I$.

\begin{figure*}
 \centering
 \vspace*{-5.cm}
    \includegraphics[width=12cm]{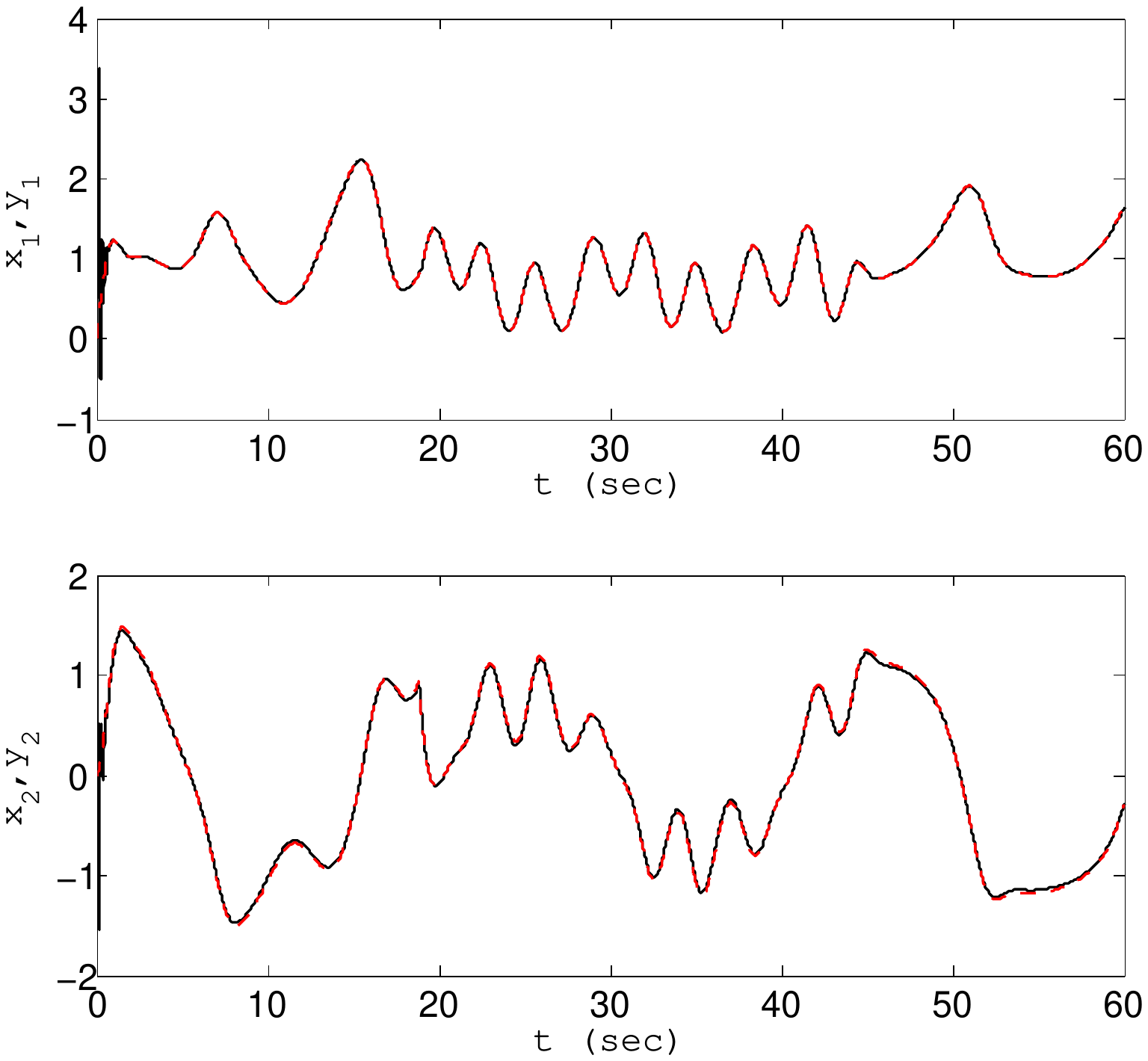}
    \vspace*{-4.cm}
    \caption{ \label{fig4} (Color online) The trajectories of states with time invariant parameters where solid
line denotes the trajectory of response system and dash line denotes the trajectory of drive system.}
 \end{figure*}

 \begin{figure}
 \centering
 \vspace*{-1.cm}
    \includegraphics[width=14cm]{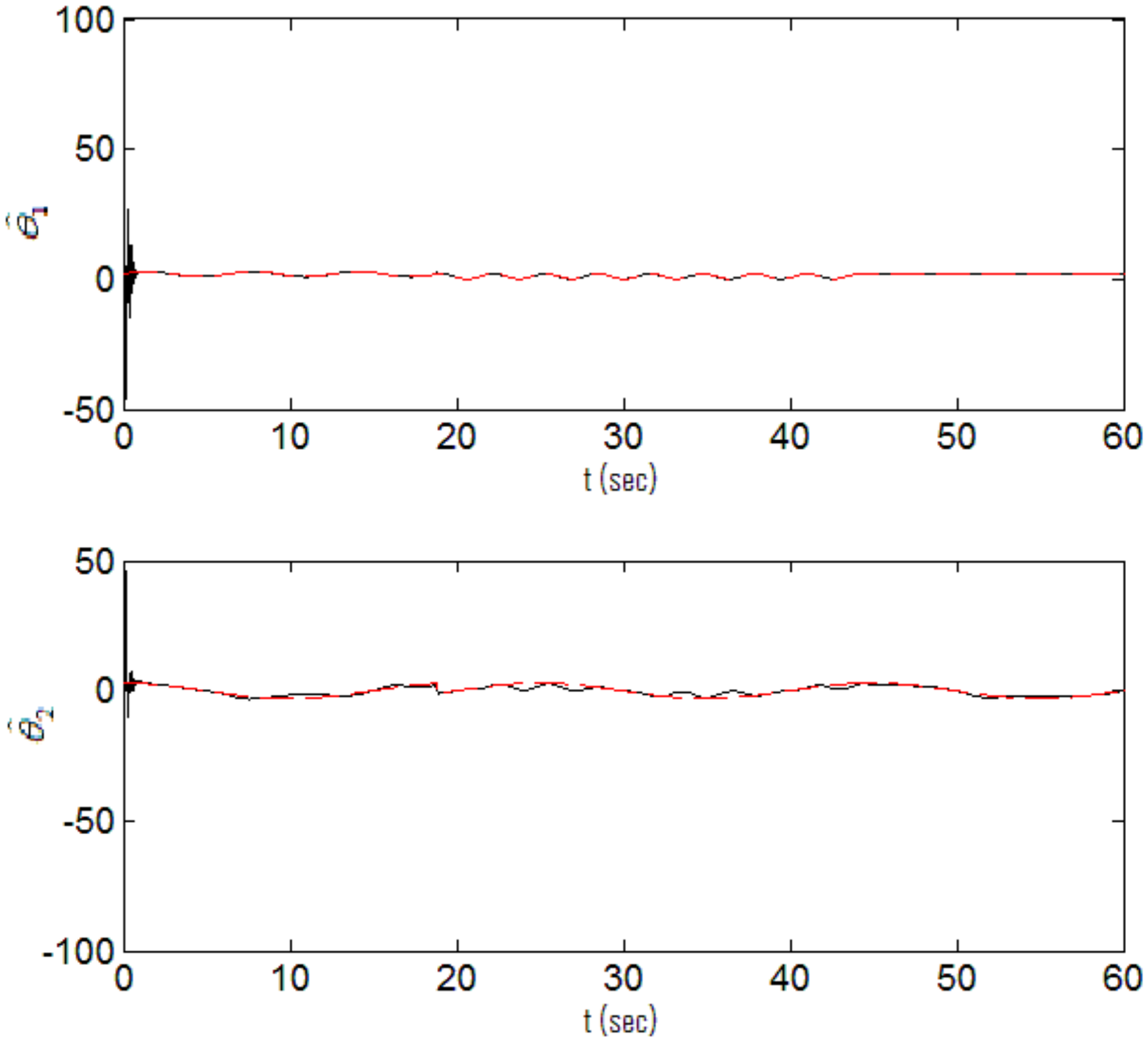}
    \vspace*{-1.cm}
    \caption{ \label{fig5} (Color online) The trajectories of estimated constant parameters where solid line
denotes the trajectory of estimated parameters and dash line denotes the trajectory of their
true values.}
 \end{figure}

 \begin{figure*}
 \centering
    \includegraphics[width=14cm]{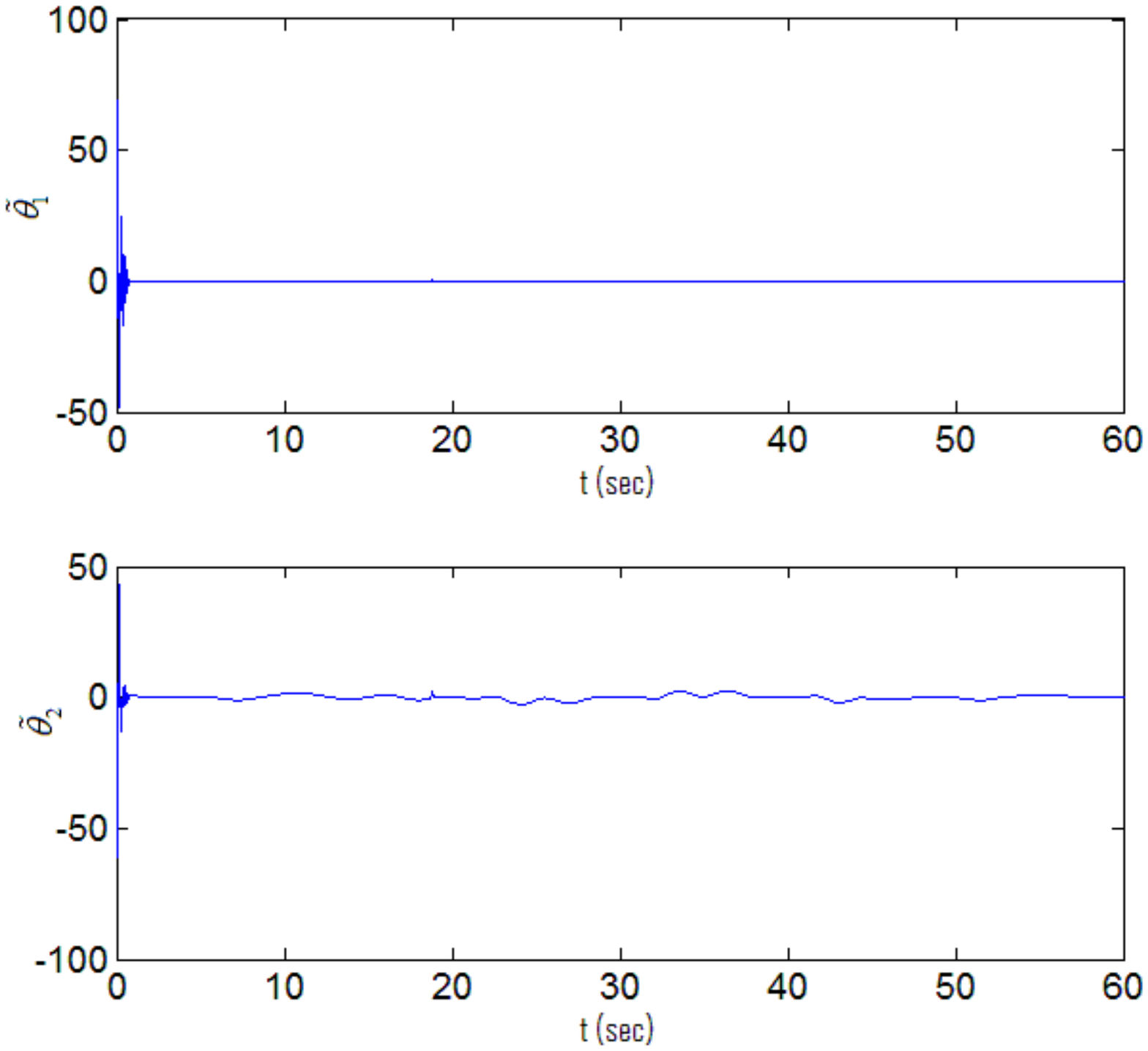}
    \vspace*{-1.cm}
    \caption{ \label{fig6} (Color online)The trajectories of estimated errors.}
 \end{figure*}

Fig.\ref{fig4} shows the trajectories of systems given in Eqs.\eqref{system2} with initial conditions provided in Eqs.\eqref{initial2}. The estimated parameters $\hat{\theta}_1$ and $\hat{\theta}_2$ are shown in Fig.\ref{fig5}, where Fig.\ref{fig6} depicts the estimated error.

It is clear from Fig.\ref{fig5} and Fig.\ref{fig6} that NRHC is able to perform as desired, and estimated time-varying parameters converge to their true values.

\section{Conclusions}\label{sec:conclusions}
In this study, we considered a real-time nonlinear receding horizon based control algorithm to estimate unknown time-varying parameters in nonlinear systems. For this purpose, we demonstrated the fact there exists a Lyapunov function which could be used to achieve asymptotic stability and provide stability guarantees for such an approach. Two benchmark examples demonstrates successful results of applicability of such concept on nonlinear systems.


\begin{thebibliography}{}

\bibitem{adetola2009} Adetola V, Dehaan D, Guay M. 2009. Adaptive model predictive control for constrained nonlinear systems. Systems and Control letters 58, 320-326.

\bibitem {ho2010} Ho W, Chou J, Guo C. 2010. Parameter identification
of chaotic systems using improved differential evolution algorithm. Nonlinear Dyn. 61, 29-41.

\bibitem {wang2012} Wang L, Xu Y. 2012. An effective hybrid biogeography-based
optimization algorithm for parameter estimation of chaotic
systems. Expert Syst. Appl. 38, 15103-15109.

\bibitem {lin2014} Lin J, Chen C. 2014. Parameter estimation of chaotic systems by an oppositional seeker optimization algorithm. Nonlinear Dyn. 76, 509-517.

\bibitem{bellsky2014} Bellsky T., Kostelich EJ, Mahalov A. 2014. Kalman filter data assimilation: Targeting observations and parameter estimation. Chaos 24, 024406.

\bibitem{ding2014} Ding F. 2014. State filtering and parameter estimation for state space systems with scarce measurements. Signal Processing 104, 369-380.

\bibitem{guay2014acc} Moshksar E, Guay M. 2014. Invariant Manifold Approach for Adaptive Estimation of the
Time-Varying Parameters for a Class of Nonlinear Systems. In Proceeding of the ACC, 2014.

\bibitem{tutsoy2015} Tutsoy O, Brown M. 2015. Chaotic dynamics and convergence analysis of temporal difference
algorithms with bang-bang control. Optim. Control Appl. Meth. DOI: 10.1002/oca.2156.

%
%
%

\bibitem{wu1993} Wu CW, Chua L. 1993. A unified framework for synchronization and control of dynamical systems. Int. J. Bifur. Chaos 3, 1619-1627.

\bibitem {li2004} Li LX, Peng HP, Wang X, Yang YX. 2004. Comment on two papers of chaotic synchronization. Phys. Lett. A 333, 269-270.
\bibitem {huang2006pre} Huang D. 2006. Adaptive-feedback control algorithm. Phys. Rev. E 73, 066204-066211.

\bibitem {li2007} Li R, Xu W, Li S. 2007. Adaptive generalized projective synchronization in different chaotic systems based on parameter identification. Phys. Lett. A 367, 199-206.

\bibitem {sun2009} Sun F, Peng HP, Luo Q, Li LX, Yang YX. 2009. Parameter identification and projective synchronization between different chaotic systems. Chaos 19, 023109.

\bibitem{sun2012} Sun Z, Si G, Min F, et al. 2012. Adaptive modified function projective
synchronization and parameter identification of uncertain hyperchaotic (chaotic)
systems with identical or non-identical structures. Nonlinear Dyn. 68, 471-486.

\bibitem{dochain2003} Dochain D. 2003. State and parameter estimation in chemical and biochemical processes: a tutorial. Journal of Process Control 13, 801-818.

\bibitem{bravo2006} Bravo JM, Alamo T, Camacho EF. 2006. Bounded error identification
of systems with time-varying parameters. IEEE Trans. Autom. Control 51, 1144-1150.

\bibitem{butt2013} Butt QR, Bhatti AI, Mufti MR, et al. 2013.
Modeling and online parameter estimation of intake manifold in
gasoline engines using sliding mode observer. Simulation Modelling
Practice and Theory 32, 138-154.

%
%
\bibitem{wang2014nd} Wang D, Ding F, Liu X. 2014. Least squares algorithm for an input nonlinear system with a dynamic subspace state space model.  Nonlinear Dyn. 75, 49-61.


\bibitem{bryson1975} Bryson AE, Ho YC. 1975. Applied optimal control. Secs. 2.3 and 6.3. Hemisphere, New York.

\bibitem{thomas1975} Thomas Y. 1975. Linear quadratic optimal estimation and control with receding horizon. Electron. Lett. 11, 19-21.

\bibitem{chen1982} Chen C, Shaw L. 1982. On receding horizon feedback control. Automatica 18, 349-352.

\bibitem{kwon1983} Kwon W, Bruckstein A, Kailath T. 1983. Stabilizing state-feedback design via the moving horizon method. Int. J. Control 37, 631-643.

\bibitem{he2008} He DF, Ji HB. 2008. Constructive model predictive control for constrained nonlinear systems. Optim. Control Appl. Meth. 29, 467-481.

\bibitem{mayne1990} Mayne D, Michalska H. 1990. Receding horizon control of nonlinear systems. IEEE Trans. Autom. Control 35, 814-824.

\bibitem{edwards1973} Edwards, C. H. 1973. Advanced Calculus of Several Variables, Theorem-3.4, pg.85, Courier Corporation

\bibitem{mayne1990b} Mayne D, Michalska H. 1990. An implementable receding horizon controller for stabilization of nonlinear systems. In proc. 29th IEEE Conf. Decision and Control,  Honolulu, HI, 3396-3397.


\bibitem {li2010pre} Peng HP, Li LX, Yang YX, Liu F. 2010. Parameter estimation of dynamical systems via a chaotic ant swarm. Phys. Rev. E 81, 016207.

\bibitem {li2006} Li LX, Yang YX, Peng HP, Wang X. 2006. Parameters identification
of chaotic systems via chaotic ant swarm. Int J Bifurcat Chaos 16, 1204-1211.

\bibitem {yu2007pre} Yu W, Chen G, Cao J, et al. 2007. Parameter identification of dynamical systems from time series. Phys. Rev. E 75, 067201.

\bibitem{parlitz1996} Parlitz U, Junge L, Kocarev L. 1996. Synchronization-based parameter estimation
from time series. Phys. Rev. E 54, 6253-6259.

\bibitem{he2007} He Q, Wang L, Liu B. 2007. Parameter estimation for chaotic
systems by particle swarm optimization. Chaos Solitons Fractals 34, 654-661.



\bibitem{zhou1996} Zhou Q, Cluett W. 1996. Recursive Identification of Time-varying
Systems via Incremental Estimation. Automatica 32, 1427-1431.


\bibitem{kabamba1987} Kabamas P, Longman R, Jian-Guo S. 1987. A homotopy approach to the feedback stablilization of linear systems. J. Guidance, Control, Dyn. 10, 422-432.

\bibitem{ohtsuka1994} Ohtsuka T, Fujii HA. 1994. Stabilized continuation method for solving optimal control problems. J. Guidance, Control, Dyn. 17, 950-957.

\bibitem{ohtsuka1997} Ohtsuka T, Fujii HA. 1997. Real-Time Optimization Algorithm for Nonlinear Receding-Horizon Control. Automatica 33, 1147-1154.

\bibitem{ohtsuka1998} Ohtsuka T. 1998. Time-variant receding-horizon control of nonlinear systems. J. Guidance, Control, Dyn. 21, 174-176.

\bibitem{lorenz1963} Lorenz E. 1963. Deterministic Nonperiodic Flow.  J. Atmos. Sci. 20, 130-141.

\end{thebibliography}

\end{document}